\documentclass[12pt,oneside]{amsart}
\usepackage{amsfonts,mathrsfs}
\usepackage{amssymb,amsmath}
\usepackage{blindtext}
\usepackage[T1]{fontenc}
\usepackage[utf8]{inputenc}
\usepackage{blkarray}

\begin{document}
 \begin{center}
{\bf\Large Matrices over maximal orders in cyclic division algebras over $\mathbb Q$  as sums of squares and cubes}\newline\\
 S. A. Katre  and Deepa Krishnamurthi \newline\\
{\bf Abstract}\\
\end{center}
It is known that every $n \times n$  matrix over the maximal order in an algebraic number field is a sum of $k$-th powers in various cases  if a discriminant condition is satisfied. It has been proved by Wadikar and Katre that for every matrix of size $\geq 2$ over maximal orders in rational quaternion division algebras  is a sum of squares and cubes. In this paper we consider  cyclic division algebras over $\mathbb Q$ of odd prime degree  and show that under some conditions every matrix of size $\ge 2$ over these noncommutative rings is a sum of squares and a  sum of cubes. 
 \\[3mm]
\\[3mm]
\newline
\section{\bf Introduction}
Let $n$ be an odd prime and $p=hn+1$. Let $g$ be a primitive root modulo $p$, $\xi$ be a primitive $p$-th root of unity. Let $\eta =\sum\limits_{r=0}^{h-1}\xi^{g^{nr}}$ be the Gaussian period and  $Z=\mathbb Q(\eta)$. The automorphism group of $\mathbb Q(\xi)$ over $\mathbb Q$ is cyclic and generated by the automorphism $\xi \rightarrow \xi^g$. Also the automorphism group of $\mathbb Q(\eta)$ over $\mathbb Q$ is cyclic generated by the automorphism $S$ which is restriction of automorphism $\xi \rightarrow \xi^g$.\\\\ {\bf Definition:} The $h$-nomial periods  are complex numbers $\eta_s=\eta_s(g)=\sum\limits_{r=0}^{h-1}exp(2 \pi i g^{kr+s}/ p),  s=0,1,2,\ldots, n-1$. \\$\eta_s=\eta^{S^s}$ is the element of $Q(\xi)$ corresponding to $\eta$ under the automorphism $S^s$. \\ $\eta_s=\sum\limits_{r=0}^{h-1}exp(2 \pi i g^{nr+s}/ p) = \sum\limits_{r=0}^{h-1}\xi^{g^{nr+s}} = \sum\limits_{r=0}^{h-1} (\xi^{g^{nr}})^{g^s}$. Hence we get $n$ automorphisms from $\mathbb Q(\eta)\rightarrow \mathbb Q(\eta)$  given by $\eta \rightarrow \eta_{s}$ induced from  $\xi\rightarrow \xi^{g^s}, s=0$ to $n-1$.\\\\ {\bf Definition:} The period polynomial of degree $n$ is given by $P_n(x)=\prod_{s=0}^{n-1}(x-\eta_s)$. It is also the minimal polynomial of $\eta$.\\\\ For a generator $g$ of $F_p^*$ the $n^2$ cyclotomic numbers of order $n$ are defined by $(i, j)_n$= cardinality of $X_{i,j}$, where  $X_{i,j}=\{v \in F_q \setminus \{0,1\}  | ind_g v\equiv i \pmod n, ind_g (v + 1) \equiv j \pmod n\}$.\\\\The following result gives us an explicit formula for the period polynomial $P_n(x)$ in terms of cyclotomic numbers of order $n$ .\\\\{\bf Dickson} {\rm\cite{lk}}: $P_n(x)=(-1)^{n-1} det(a_{ij})$, \\$ a_{11}= 1+x,\\ a_{1j}=1, ~~2 \le  j \le n\\ a_{i1}=h\theta_i+(i-1,0)_n x; ~~2 \le i \le n\\ a_{ii}=(i-1,i-1)_n -x;~~2 \le i \le n\\ a_{ij}=(i-1, j-1)_n; ~~ 2 \le i \le n, 2 \le j \le n, i \neq  j,$ where $\theta_i =1$ if $h$ is even  or $i=0$ and $h$ is odd and $\theta_i =0$ otherwise.\\\\ {\bf Example:}  For $n=3$ \\ $P_3(x)= (-1)^2 $ det $ (a_{ij})_ {3\times3}$\\ $a_{11}=1+x\\ a_{12}=1\\ a_{13}=1\\ a_{21}=(1,0)_3 x\\ a_{22}=(1,1)_3 -x\\ a_{23}=(1,2)_3\\ a_{31}=(2,0)_3 x\\ a_{32}=(2,1)_3\\ a_{33}=(2,2)_3-x$\\\[ P_3(x)=
 \begin{vmatrix}
1+x & 1 & 1\\
(1,0)_3x & (1,1)_3-x & (1,2)_3\\
(2,0)_3x & (2,1)_3 & (2,2)_3-x
\end{vmatrix}.\]\\ Now cyclotomic matrix of order 3 is 
\[
\begin{blockarray}{ccc}
\begin{block}{(ccc)}
   A & B & C  \\
   B & C & D  \\
   C & D & B  \\
\end{block}
\end{blockarray}
 \] whose $(i,j)^{th}$ entry is given by $(i, j)_3$\\ where,\\ $9A= p+r_3 - 8\\ 18B= 2p - 4 - r_3 -9d\\ 18C= 2p-4 - r_3 +9d\\ 9D=p+1+r_3$\\ where $r_3$ is uniquely determined by $4p= r_3^2 +27 d^2, r_3\equiv 1\pmod 3$ and $d$ is determined upto sign and sign depends on the generator of $F_p^*$. Here $r_3=1$, we take the generator $g$ for which $d=1$.\\\\ Using above relations  for $p=7, n=3$  we have $d=\pm1, r_3=1$, we take $d=1$. Also $18B=0, B =(1,0)_3 $ hence $B=(1,0)_3=0$ and $18C=18$, so $C=1$, $(1,1)_3=C, (2,2)_3 =B=0,~ 9D=9$, hence $D=1$. So $D=(1,2)_3=1$.\\\[ P_3(x)=
 \begin{vmatrix}
1+x & 1 & 1\\
0     & 1-x & 1\\
x &  1 &  -x
\end{vmatrix}.\]\\ $P_3(x)= x^3+x^2-2x-1$.\\\\ {\bf Myerson {\rm\cite{my}}}: In general for $n=3$ and any prime $p$, \\$P_3(x)= x^3+x^2 - \frac{(p-1)x}{3} - \frac{(r_3+3)p-1}{27}$.\\ The constant term of the period polynomial is the product $-\eta\cdot\eta^{S}\cdot\eta^{S^2}$, which can be expressed in terms of the cyclotomic numbers of order 3.\\ \\{\bf Matrices over Cyclic Algebras as sums of squares }\\\\{\bf Definition:} A cyclic algebra  $A$ over $\mathbb Q$ of degree $n$ and order $n^2$  has a $\mathbb Q$-basis of the form \begin{enumerate}\item $u^iz_k, i=0,1,2,\ldots, n-1 ; k=1,2,\ldots n$ \\where $z_k$ forms a $\mathbb Q$-basis of the cyclic subfield $Z$ of $A$ of degree $n$ over $\mathbb Q$  and $1, u, u^2, \ldots, u^{n-1}$ form a $Z$ basis of $A$, with the relations \item $zu=uz^S$ \\where $S$ is the generating element of the Galois group of $Z$ over $\mathbb Q$ and $z^S$ is the element of $Z$ coresponding to $z$ under the automorphism $S$, and\item $u^n=\sigma \neq 0$ in $\mathbb Q$.\\ This is called cyclic generation of $A$ and is denoted by $A=(\sigma, Z, S)$ \end{enumerate}{\bf Proposition 1 (R. S. Pierce}{\rm(\cite{p})}): For a cyclic extension $L/K$ of
degree $n$ with Galois group $Gal(L/K)$ = $<S>$. If $\alpha$, $\alpha^{2}$, $\ldots$, $\alpha^{n-1} \in K^{*}$ are not norm of some element of $L$ then $(L/K, S, \alpha)$ is a cyclic division algebra.\\\\{\bf Examples of cyclic division algebras}:\begin{enumerate}\item Consider $\mathbb Q+\mathbb Q(i)+\mathbb Q(j)+\mathbb Q(k)= \mathbb Q(i)+ j \mathbb Q(i)$, then $\mathbb Q(i)$ has degree 2 over $\mathbb Q$ and $S$  is the automorphism $i \rightarrow -i$, also $\alpha=-1$ is not norm of any element in $\mathbb Q (i)$, so by Proposition 1, $\mathbb Q$+$\mathbb Q(i)$+$\mathbb Q(j)$+$\mathbb Q(k)$ is a cyclic division algebra of degree 2 and order 4.\item Consider $\mathbb Q(\xi_7+\xi_7^{-1})$+ $u\mathbb Q(\xi_7+\xi_7^{-1})$+ $u^2 \mathbb Q(\xi_7+\xi_7^{-1})$ where $\xi_7$  is a primitive seventh root of unity and $u^3=2$.  Here  $\mathbb Q(\xi_7+\xi_7^{-1})$ has degree 3 over $\mathbb Q$. The  minimal polynomial satisfied by $(\xi_7+\xi_7^{-1})$ is $x^3+x^2-2x-1$ and $S$ is the automorphism $(\xi_7+\xi_7^{-1})\rightarrow (\xi_7^2+\xi_7^5)$. Also $\alpha =2$ and $\alpha^{2}=4$ are not norm of any element in $\mathbb Q(\xi_7+\xi_7^{-1})$. Hence by Proposition 1, $\mathbb Q(\xi_7+\xi_7^{-1})$+ $u\mathbb Q(\xi_7+\xi_7^{-1})$+ $u^2 \mathbb Q(\xi_7+\xi_7^{-1})$ is a cyclic division algebra of degree 3 and order 9.\end{enumerate}{\bf Definition:} An order in a rational semi-simple algebra $B$ is a subset $I$ of elements of $B$ with the following properties:\begin{enumerate}\item $I$ contains unity element. \item The sum difference and product of any two elements in $I$ is again in $I$. \item If $b$ is any element of $B$ there exists a rational integer $\mu $ such that $\mu b \in I$.\item There is a finite set of elements $a_1, a_2, \ldots, a_r$ such that every element $a \in I$ can be expressed as $a=n_1a_1+\cdots+n_ra_r$ where $n_1, n_2,\ldots,n_r$ are rational integers.\end{enumerate}  An order is maximal if it cannot be embedded in any other order.\\\\{\bf Theorem 1 ((Ralph Hull {\rm\cite{skb})}, Existence of canonical generations)}\\  Let $A$ be a cyclic division algebra of odd prime degree $n$ over $\mathbb Q$, let $q_1, q_2, \ldots, q_s (s\geq 2)$ be distinct rational primes at which $A$ does not split. Let $\sigma=q_1\cdot q_2 \cdots q_s$. Then there exists infinitely many rational primes $p$ with the following properties: $p\equiv1\pmod n$ and is prime to $\sigma$; $q_1, q_2, \ldots q_s$ are $n$-ic nonresidues modulo $p$ and $\sigma $ is an $n$-ic residue modulo $p$; the unique cyclic field $Z$ of degree $n$ over $\mathbb Q$, of conductor $p$, discriminant $p^{n-1}$, has an automorphism $S$ such that $(\sigma, Z, S)$ is a cyclic generation of $A$.\\\\ Such cyclic generations as in Theorem 1 are called as canonical generations.\\\\{\bf Order associated with a canonical generation:} Let $(\sigma, Z, S)$ be a cyclic generation as in Theorem 1, then the totality of all linear combinations of the basis elements $u^i z_k$ with rational integral coefficients is called order $I$ associated with the canonical generation.\\\\{\bf Definition:} The set \{$z_1, z_2,\ldots, z_n$\} of $n$ conjugates obtained by applying automorphisms of $\mathbb Q(\eta)$ is called normal basis for integers of $Z$.\\\\The following theorem gives description of maximal orders in $A$ containing $I$, see \cite{skb}. By Theorem 1 since $\sigma $ is a $n$-ic residue modulo $p$, the congruence $\lambda^ n \equiv \sigma \pmod p$ has a rational integral solution. Let $\lambda$ be such a fixed solution, let $\alpha= \beta ^{n-1}$, where $\beta=\prod_{r=0}^{h-1}(1-\xi^{g^{nr})}, ~g$ is primitive root modulo $p$, $\xi^p=1$. Let $Z=\mathbb Q(\eta)$ where $\eta=\sum\limits_{r=0}^{h-1}\xi^{g^{nr}}$ is the Gaussian period  and $y$ be defined by the equation $py=(\lambda- u) \alpha$. \\\\ {\bf Theorem 2 (Ralph Hull {\rm(\cite{skb})}, Existence of maximal orders containing I):} For a fixed integral solution $\lambda$ of $\lambda^n \equiv \sigma$ (mod $p$), the set $I(\lambda)$ of linear combinations with rational integral coefficients of the quantities $y^i z_k, i=0,\ldots, n-1; j=1,2,\ldots n$ form a normal basis of integers of $Z$, is a maximal order in $A$ containing $I$.\\\\{\bf Notation}: $ T_{2,n}=T_2$ is set of those elements of cyclic algebra which can be expressed as sums of traces of squares of $n\times n$ matrices over the cyclic algebra.\\$\mathfrak O$ is the ring of algebraic integers in $\mathbb Q(\eta)$ generated by $\eta, \eta^S, \ldots, \eta^{S^{n-1}}$\\\\As in Theorem 1, in the following theorem we consider the case where all $q_i$'s are odd. This means the cyclic algebra splits at 2. In this case $\sigma \equiv 1 \pmod 2$.\\\\{\bf Theorem 3:} Let $A$ be a rational cyclic division algebra of degree 3 and order 9 with canonical generation $(\sigma, Z, S)$ as in Theorem 1 and let $\mathfrak{m}$ be a maximal order associated with this canonical generation. If $\sigma \equiv 1 \pmod 2$ then every matrix $C \in M_n (\mathfrak{m})$ is sum of squares of matrices in  $M_n (\mathfrak{m})$.\vspace{0.5cm}\\{\bf Proof:} $\mathfrak{m}$ is a free $Z$ module of rank 9 and there are nine $Z$-basis elements of $\mathfrak{m}$ namely \begin{center}$\eta,~\eta^S,~\eta^{S^2}, ~y\cdot \eta, ~y\cdot \eta^S, ~y \cdot \eta^{S^2}, ~y^2\cdot \eta, ~y^2 \cdot \eta^S,~y^2 \cdot \eta^{S^2}$\end{center} so it suffices to prove that these basis elements are in $T_2$.\\Now $\eta=\sum\limits_{r=0}^{h-1}\xi^{g^{nr}}$, if $i=2k$ then $\xi^i=\xi^{2k}$ is already a square, hence in $T_2$. If $i=2k+1$ then $\xi^i =\xi^{2k+1}=(\xi^2)^{(\frac{2k+1+p}{2})}$, hence a square, so in $T_2$. Hence $\eta \in T_2$. Since $\eta, ~\eta^S, ~\eta ^{S^2}$ are conjugates, so $\eta,~\eta^S,~\eta^{S^2} \in T_2$. Next to prove $y \in T_2$. Now $py=(\lambda - u) \alpha$ and $p$ being odd prime, $p \equiv 1\pmod 2.$ Hence $y \equiv (\lambda - u) \alpha \pmod 2\equiv (\lambda - u)(a_1 \eta+a_2 \eta^S+a_3 \eta^{S^2})\pmod 2,~ a_i \in \mathbb{Z}$. Now $u^3 \equiv \sigma \equiv 1\pmod 2$, hence $u^4 \equiv u\pmod 2$. So $u \in T_2$. Hence it suffices to prove that $u \cdot \eta$, $u\cdot \eta^S$ and $u \cdot \eta^{S^2} \in T_2$. Now \begin{center}$(u+\eta^S)^2= u^2+u\cdot \eta^S+\eta^S\cdot u+ (\eta^S)^2$\\=$ u^2+u\cdot \eta^S+u \cdot \eta^{S^2}+(\eta^S)^2 $\\$=u^2+u \cdot (\eta^S+ \eta^{S^2})+(\eta^S)^2$\\ =$ u^2+u \cdot (-1- \eta)+(\eta^S)^2$.\end{center} Hence $u \cdot \eta \in T_2$. Similarly by expanding $(u +\eta^{S^2})^2$ and $(u+\eta)^2$ we get $u \cdot \eta^S$ and $u \cdot \eta^{S^2} \in T_2$. Hence $y \in T_2$. Now to prove that $y\cdot \eta,~ y \cdot \eta^S $ and $y \cdot \eta^{S^2} \in T_2$. Consider $y\cdot \eta=(\lambda- u) \alpha\cdot\eta=(\lambda- u)(b_1\cdot \eta+b_2\cdot \eta^S+b_3\cdot \eta^{S^2}),~b_i \in \mathbb{Z} $, so it suffices to prove that  $u\cdot \eta, u \cdot \eta^S$ and $u \cdot \eta^{S^2} \in T_2$ which is already proved. Next to prove $y^2\cdot \eta, y^2 \cdot \eta^S,  y^2 \cdot \eta^{S^2} \in T_2$. For this $y^2 \cdot \eta = [(\lambda- u)(a_1\cdot \eta+a_2\cdot \eta^S+a_3\cdot \eta^{S^2})]\times[(\lambda- u)(a_1\cdot \eta+a_2\cdot \eta^S+a_3\cdot \eta^{S^2})]\cdot\eta$ The product on R.H.S is a $ \mathbb Z$  linear combination of terms of the type $\eta^{S^i}\cdot \eta^{S^j}\cdot\eta, ~u\cdot\eta^{S^i}\cdot\eta^{S^j}\cdot \eta, ~\eta^{S^i}\cdot u\cdot\eta^{S^j}\cdot\eta$ and $u \cdot\eta^{S^i}\cdot u\cdot\eta^{S^j}, 0\leq i, j \leq 2 $.\begin{center} Now $\eta^{S^i}\cdot u= u\cdot \eta^{S^{i+1}}$\end{center} We see that $y^2 \cdot \eta = c_1+ u \cdot c_2 + u^2\cdot c_3,  c_1, c_2, c_3 \in \mathfrak O$, hence it suffices to prove that $u^2 \cdot \eta \in T_2$.\begin{center} Now $(u^2+\eta^{S^2})^2= u^4+u^2 \cdot \eta^{S^2}+\eta^{S^2}\cdot u^2+(\eta^{S^2})^2$\\$= u^4+ u^2\cdot(-1-\eta)+(\eta^{S^2})^2$\end{center} Hence $u^2 \cdot \eta \in T_2$. Thus $y^2\cdot \eta \in T_2 $. By similar argument $ y^2 \cdot \eta^S,~ y^2 \cdot \eta^{S^2} \in T_2$. Hence all basis elements of $\mathfrak M$ are in $T_2$.\qed\\\\ We now consider condition for general odd prime $n$. Our proof requires that $\eta^{S^{n-1}}\cdot \eta^{S^{n-2}}\cdots \eta \equiv 1\pmod 2$ i.e. constant term in the minimal polynomial of $\eta$ is an odd number.\\\\{\bf Theorem 4:}  Let $A$ be a rational cyclic division algebra of degree $n$ and order $n^2, ~n$  is an odd prime  with canonical generation $(\sigma, Z, S)$ as in Theorem 1 and let $\mathfrak{m}$ be a maximal order associated with this canonical generation satisfying the following conditions \begin{enumerate}\item $\sigma \equiv 1\pmod 2$ \item $N_{\mathbb{Q}(\eta)/\mathbb{Q}}(\eta) \equiv 1\pmod 2.$\end{enumerate} Then every matrix $C \in M_n (\mathfrak{m})$ is a sum of squares of matrices in $M_n (\mathfrak{m}).$\\\\{\bf Proof:}$~\mathfrak {m}$ is a free $Z$ module of rank $n^2$ and there are $n^2$ basis elements of $\mathfrak{m}$ namely \begin{center}$\eta,~\eta^S,~\ldots, ~\eta^{S^{n-1}},~y\cdot \eta,~y\cdot n^S, \ldots,~y \cdot \eta^{S^{n-1}}, ~y^2\cdot \eta, ~y^2 \cdot \eta^S,\ldots, \newline y^2 \cdot \eta^{S^{n-1}},~\ldots, ~y^{n-1} \cdot \eta,~ y^{n-1}\cdot \eta^S, ~\ldots,~y^{n-1} \cdot \eta^{S^{n-1}}$\end{center} so it suffices to prove that these basis elements are in $T_2$.\\ Now $\eta = \sum\limits_{r=0}^{h-1}\xi^{g^{nr}}$\\ {\bf Case 1}: If $i=2k$ then $\xi^i=\xi^{2k}$ it is already a square, hence  $\eta \in T_2$.\\ {\bf Case 2}: If $i=2k+1$ then $\xi^i =\xi^{2k+1}=(\xi^2)^{(\frac{2k+1+p}{2})},$ hence a square, so $\eta \in T_2$. \\Since $\eta, ~\eta^S,~ \eta ^{S^2}, ~\ldots,~ \eta^{S^{n-1}}$ are conjugates, so $\eta, ~\eta^S, ~\eta^{S^2}~\ldots,~\eta^{S^{n-1}} \in T_2$. Next to prove $y \in T_2$. Now $py=(\lambda - u) \alpha$ and $p$ being odd prime, $p \equiv 1\pmod 2.$ Hence $y \equiv (\lambda - u)\alpha$(mod 2)$ \equiv (\lambda - u)(a_1 \eta+a_2 \eta^S+a_3 \eta^{S^2}+\cdots+a_n \eta^{S^{n-1}})\pmod 2, ~a_i \in \mathbb{Z}.$ For this we first prove $u \in T_2.$ Now $u^n \equiv \sigma \equiv 1 \pmod 2$, hence $u^{n+1} \equiv u\pmod 2.$ Now $n$ being  odd, $n+1$  is even number. So $u \in T_2$. \\Next to prove that $u\cdot \eta,~ u\cdot \eta^S,~\ldots,~ u \cdot \eta^{S^{n-1}} \in T_2,$ it suffices to prove that $u \cdot \eta \in T_2$.  Now $(u\cdot\eta)^{n} \equiv u^n\cdot\eta^{S^{n-1}}\cdot \eta^{S^{n-2}}\cdots\eta \equiv \sigma\cdot \eta^{S^{n-1}}\cdot \eta^{S^{n-2}}\cdots\eta \equiv 1.1 \equiv 1 \pmod 2$. So $(u \cdot \eta)^{n+1} \equiv u\cdot \eta\pmod 2.$ Now $n$ is odd, so $n+1$ being even number, $u \cdot \eta \in T_2$. By similar argument we can prove that $u \cdot \eta^S$ and $u \cdot \eta^{S^2}, \ldots, u\cdot \eta^{S^{n-1}} \in T_2$. Hence $y \in T_2$.\newline Next to prove $y\cdot\eta, y\cdot \eta^S, \ldots, y\cdot \eta^{S^{n-1}} \in T_2$. \\Now $y\cdot \eta=(\lambda- u) \alpha \cdot \eta = (\lambda- u)(b_1 \eta+b_2 \eta^S+\cdots+b_n \eta^{S^{n-1}}) $ where $b_i \in \mathbb{Z}$, so suffices to prove that $u\eta, \ldots, u\eta^{S^{n-1}}\in T_2$ which is already proved.\newline Next to prove  $ y^2\cdot \eta,~y^2 \cdot \eta^S,\ldots, y^2 \cdot \eta^{S^{n-1}}\in T_2$. Now $y^2 \cdot \eta$ \begin{center}$= [(\lambda- u)(a_1\cdot \eta+a_2\cdot \eta^S+a_3\cdot \eta^{S^2}+\cdots+a_n \eta^{S^{n-1}})]\times [(\lambda- u)(a_1\cdot \eta+a_2\cdot \eta^S+a_3\cdot \eta^{S^2}+\cdots+\eta^{S^{n-1}})] \cdot \eta$ \end{center} The product on R.H.S is a $\mathbb Z$- linear combination of terms of the type \begin{center}$\eta^{S^i}\cdot \eta^{S^j}\cdot\eta, ~u\cdot\eta^{S^i}\cdot\eta^{S^j}\cdot \eta,~ \eta^{S^i}\cdot u\cdot\eta^{S^j}\cdot\eta $ and $u \cdot\eta^{S^i}\cdot u\cdot\eta^{S^j} , 0\leq i,~j \leq n-1 $.\end{center} Now \begin{center}$\eta^{S^i}\cdot u= u\cdot \eta^{S^{i+1}}$ \end{center} We see that $y^2 \cdot \eta = c_1+ u \cdot c_2 + u^2\cdot c_3 , c_1, c_2, c_3, \in \mathfrak O$ so it suffices to show that $u^2\cdot\eta \in T_2$ .\begin{center} Consider $(u^2\cdot \eta)^n\equiv u^{2n}\cdot\eta^{S^{2(n-1)}} \cdot \eta^{S^{2(n-2)}}\cdots\eta \equiv 1.1 \equiv 1\pmod 2.$\end{center} Hence $(u^2\cdot \eta)^{n+1}\equiv u^2 \cdot \eta\pmod 2.$ Since $n+1$ is even, this implies $u^2\cdot \eta \in T_2$.  Now to prove that  for $3 \leq  i  \leq  n-1,~ y^i \cdot \eta \in T_2$, it suffices to prove that $u^i \cdot \eta, u^i\cdot \eta^S, \ldots, u^i \cdot \eta^{S^{n-1}} \in T_2$ (for clarity purpose cases $i=0, 1, 2$ have been separately done before). Since $n$ is odd and $\eta \cdot u^i= u^i \cdot \eta^{S^i}$ we have $ (u^i \cdot \eta)^n\equiv u^{in}\eta^{S^{i(n-1)}}\cdot \eta^{S^{i(n-2)}}\cdots \eta^{S^{2i}}\cdot \eta^{S^i} \cdot \eta$. As $S$ is of prime order $n$ and generator of the automorphism group  of $\mathbb Q(\eta)$ over $\mathbb Q$, so $S^i$ is also a generator of $\mathbb Q(\eta)$ over $\mathbb Q$, hence $\eta^{S^{i(n-1)}},\eta^{S^{i(n-2)}},\ldots, \eta^{S^{2i}}, \eta^{S^i},\eta$ are  $\eta, \eta^S, \cdots, \eta^{S^{n-1}}$ in some order. \\So we have, \begin{center} $ (u^i \cdot \eta)^n\equiv u^{in}\cdot\eta^{S^{i(n-1)}}\cdot \eta^{S^{i(n-2)}}~\cdots~ \eta^{S^{2i}}\cdot \eta^{S^i }\cdot \eta$ \\ $\equiv u^{in}\cdot \eta \cdot \eta^S\cdots\eta^{S^{n-1}}\pmod 2$\\$\equiv \sigma^i\cdot \eta \cdot \eta^S\cdots\eta^{S^{n-1}}\pmod 2$\\$\equiv 1.1\equiv 1\pmod 2.$\end{center} Hence $(u^i\cdot \eta)^{n+1}\equiv u^i\cdot\eta\pmod 2.$ Since $n$ is odd, so $u^i \cdot \eta \in T_2$. Thus $y^i \cdot \eta \in T_2$. So all basis elements of $\mathfrak M \in T_2$.\\\\{\bf Matrices over cyclic division algebras as sums of cubes}\\\\{\bf Notation}: $ T_{3,n}=T_3$ be the set of those elements of the cyclic algebra $A$ which can be written as sums of traces of cubes of $n\times n$ matrices over the cyclic algebra. \\\\{\bf Theorem 5:} Let $A$ be a rational cyclic division algebra of degree $n$ and order $n^2, ~n$  is an odd prime  with canonical generation $(\sigma, Z, S)$ as in Theorem 1 and let $\mathfrak{m}$ be a maximal order associated with this canonical generation satisfying following conditions \begin{enumerate}\item $3 \nmid \sigma$  and \item $3\nmid N_{\mathbb{Q}(\eta)/\mathbb{Q}}(\eta).$ \end{enumerate}Then every matrix $C \in M_n (\mathfrak{m})$ is a sum of cubes of matrices in  $M_n (\mathfrak{m})$\\\\{\bf Proof:} $\mathfrak {m}$ is a free $Z$ module of rank $n^2$ and there are $n^2$ basis elements of $\mathfrak{m}$ namely \begin{center}$\eta, ~\eta^S,~\ldots, ~\eta^{S^{n-1}}, ~y\cdot \eta,~ y\cdot n^S,~ \ldots, ~y \cdot \eta^{S^{n-1}}, ~y^2\cdot \eta, ~y^2 \cdot \eta^S,~\ldots,\newline y^2 \cdot \eta^{S^{n-1}},~\ldots, y^{n-1} \cdot \eta, ~y^{n-1}\cdot \eta^S, \ldots,y^{n-1} \cdot \eta^{S^{n-1}} $\end{center}  so we must  prove that these basis elements are in $T_3$. \\Now $\eta=\sum\limits_{r=0}^{h-1}\xi^{g^{nr}}, \xi^p=1$, so $\xi^{p+1}=\xi$. First we  prove $\eta \in T_3$.  Now  consider the congruence $3x\equiv (p+1)\pmod p.$ This always has a solution as $(3,p)=1$. So $\xi^{p+1}=\xi^{3x}$ . Hence $\xi \in T_3$ . Also $(\xi^i)^{p+1}=(\xi^i)^{3x}$. \\Thus $\eta \in T_3$. Since $\eta, ~\eta^S, ~\eta ^{S^2}, ~\ldots, ~\eta^{S^{n-1}}$ are conjugates, hence $\eta, ~\eta^S,~ \eta^{S^2}, \ldots, ~\eta^{S^{n-1}} \in T_3$. \\Next to prove $y \in T_3$. Now $py=(\lambda - u) \alpha$ and $p$ being odd prime, $ p > 3$, $p \equiv 1\pmod 3$ or $p \equiv -1 \pmod 3$ . If $p \equiv 1\pmod 3$ then $y=(\lambda - u) \alpha=(\lambda - u)(a_1 \eta+a_2 \eta^S+a_3 \eta^{S^2}+\cdots+a_n \eta^{S^{n-1}})$ and if $p\equiv-1\pmod 3$ then $- y=(\lambda - u) \alpha=(\lambda - u)(a_1 \eta+a_2 \eta^S+a_3 \eta^{S^2}+\cdots+a_n \eta^{S^{n-1}})$. So it suffices to prove that  $u\cdot \eta^{S^i}, i=0, 1, \ldots, n-1 \in T_3$. First we prove $u \in T_3$. We consider following cases: \\{\bf Case 1}: If $n=3k+2$ then $u^{3k+2}\equiv u^{n}\equiv \sigma\equiv 1\pmod 3,$ so $u^{n+1}=u^{3k+3}\equiv u\pmod 3.$ Hence $u \in T_3.$\\{\bf Case 2}: If $n=3k+1$, then $u^{n+1}\equiv u\equiv u^{3k+2}\pmod 3,$ so $u\equiv (u^{3k+2})^{3k+2}\equiv u^{3k'+1}\pmod 3.$ Now $u^n \equiv \sigma \equiv 1 \pmod 3,$ hence $ u^{3(k'-k)}\equiv u^{(3k'+1) - (3k+1)}\cdot u^{3k+1}\equiv u ^{3k'+1}\equiv u\pmod 3.$ So $u \in T_3$. \\Next to prove $y \in T_3$. \\Now $py=(\lambda - u) \alpha$ and $p \equiv 1\pmod 3$ or $p \equiv -1\pmod 3$. If $p \equiv 1\pmod 3$, then $y=(\lambda - u) \alpha=(\lambda - u) (b_1\eta + b_2\eta^S+\cdots+b_n\eta^{S^{n-1}})$ or if $p \equiv -1\pmod 3$,  $- y =(\lambda - u) \alpha= (\lambda - u)(b_1\eta + b_2\eta^S+\cdots+b_n\eta^{S^{n-1}}),$ so suffices to prove that $ u\cdot \eta \in T_3$. Now $(u\cdot\eta)^{n} \equiv u^n\cdot\eta^{S^{n-1}} \cdot \eta^{S^{n-2}}\cdots\eta \equiv 1.1 \equiv 1\pmod 3.$ \newline{\bf Case 1}: If $n=3k+2$, so $(u\cdot\eta)^{n+1} \equiv u \cdot \eta \equiv (u\cdot \eta)^{3(k+1)}$(mod 3), so $u\cdot \eta \in T_3$. \newline{\bf Case 2}: If $n=3k+1$, then $(u\cdot \eta)^{n+1}\equiv u\eta\equiv (u\cdot\eta)^{3k+2}\pmod 3$, so $u\cdot \eta\equiv ((u\cdot \eta)^{3k+2})^{3k+2}\equiv (u\cdot \eta)^{3k'+1}\pmod 3$. Now $(u\cdot \eta)^n \equiv \sigma\cdot \eta\cdot\eta^S~\cdots~\eta^{S^{n-1}} \equiv 1\pmod 3$, hence $(u\cdot \eta)^{3(k'-k)}\equiv (u\cdot \eta)^{(3k'+1) - (3k+1)} \cdot ( u \cdot \eta)^{3k+1}\equiv (u\cdot \eta)^{3k'+1} \equiv u\cdot \eta \pmod 3$. So $u \cdot \eta \in T_3$. \\ By similar argument replacing $u\cdot\eta$ by $u\cdot \eta^{S^{i}}, 1 \leq i \leq n-1$, we get $u\cdot \eta^{S^{i}} \in T_3$, hence $y \in T_3$ . Next to prove $y\cdot\eta, y\cdot \eta^S, \ldots, y\cdot \eta^{S^{n-1}} \in T_3$. Now $y\cdot \eta=(\lambda- u) \alpha \cdot \eta = b_1 u\cdot \eta+b_2 u \cdot\eta^S+\cdots+b_n u\cdot \eta^{S^{n-1}}$ where $b_i \in \mathbb{Z}$, so suffices to prove that $u\cdot \eta, \ldots, u\cdot \eta^{S^{n-1}}\in T_3$ which is already proved.  Next to prove  $ y^2\cdot \eta, y^2 \cdot \eta^S,\ldots,  y^2 \cdot \eta^{S^{n-1}}\in T_3$. Now $y^2 \cdot \eta$ \begin{center}$= [(\lambda- u)(a_1\cdot \eta+a_2\cdot \eta^S+a_3\cdot \eta^{S^2}+\cdots+a_n \eta^{S^{n-1}})]\times [(\lambda- u)(a_1\cdot \eta+a_2\cdot \eta^S+a_3\cdot \eta^{S^2}+\cdots+a_n\eta^{S^{n-1}})] \cdot \eta$ \end{center} The product on R.H.S is a $\mathbb {Z}$- linear combination of terms of the type \begin{center}$\eta^{S^i}\cdot \eta^{S^j}\cdot\eta, ~u\cdot\eta^{S^i}\cdot\eta^{S^j}\cdot \eta,~ \eta^{S^i}\cdot u\cdot\eta^{S^j}\cdot \eta $ and $u \cdot\eta^{S^i}\cdot u\cdot\eta^{S^j}, 0\leq i,j \leq n-1 $.\end{center} Now \begin{center}$\eta^{S^i}\cdot u= u\cdot \eta^{S^{i+1}}$\end{center} so it suffices to show that $u^2\cdot\eta,~ u^2\cdot\eta^S, \ldots, u^2\cdot\eta^{S^{n-1}} \in T_3.$ Now $(u^2\cdot\eta)^{n} \equiv u^{2n}\cdot\eta^{S^{2(n-1)}} \cdot \eta^{S^{2(n-2)}}\cdot\eta \equiv 1.1 \equiv 1 \pmod 3.$\newline{\bf Case 1}: If $n=3k+2$, so $(u^2\cdot\eta)^{n+1} \equiv u^2 \cdot \eta \equiv (u^2\cdot \eta)^{3(k+1)}\pmod 3$, so $u^2\cdot \eta \in T_3$. \newline {\bf Case 2}: If $n=3k+1$, then $(u^2\cdot \eta)^{n+1}\equiv u^2\cdot \eta\equiv (u^2\cdot\eta)^{3k+2}\pmod 3,$ so $u^2\cdot \eta\equiv ((u^2\cdot \eta)^{3k+2})^{3k+2}\equiv (u^2\cdot \eta)^{3k'+1}\pmod 3.$ Now $(u^2\cdot \eta)^n  \equiv 1\pmod 3,$ hence $u^2\cdot \eta\equiv(u^2\cdot \eta)^{(3k'+1) - (3k+1)}\equiv (u^2\cdot \eta)^{3(k'-k)}\pmod 3.$ Hence $u^2 \cdot \eta \in T_3$.\\ By induction we must prove that $y^{n-1} \cdot \eta,~ y^{n-1} \cdot \eta^S, \ldots , ~y^{n-1} \cdot \eta^{S^{n-1}} \in T_3$. As above it suffices to prove that $u^{n-1}\cdot\eta \in T_3$.  \newline{\bf Case 1}: If $n=3k+2$, so $(u^{n-1}\cdot\eta)^{n+1} \equiv u^{n-1} \cdot \eta \equiv (u^{n-1}\cdot \eta)^{3(k+1)}\pmod 3,$ so $u^{n-1}\cdot \eta \in T_3$. \newline{\bf Case 2}: If $n=3k+1$, then $u^{n-1}\cdot \eta= (u^{n-1}\cdot\eta)^{3k+2}$, so \\$u^{n-1}\cdot \eta\equiv ((u^{n-1}\cdot \eta)^{3k+2})^{3k+2}\equiv (u^{n-1}\cdot \eta)^{3k'+1}\pmod 3.$ Now $(u^{n-1}\cdot \eta)^n \equiv \sigma^{n-1}\cdot\eta^{S^{n-1}}\cdot\eta^{S^{n-2}}~\cdots~\eta \equiv 1 \pmod 3,$ hence $u^{n-1}\cdot \eta \equiv (u^{n-1}\cdot \eta)^{3k'+1}\equiv (u^{n-1}\cdot \eta)^{(3k'+1) - (3k+1)}\cdot (u^{n-1}\cdot \eta)^{3k+1}\equiv (u^{n-1}\cdot \eta)^{3(k'-k)}\pmod 3.$ So $u^{n-1}\cdot \eta \in T_3$. \\Thus all basis elements of $\mathfrak{m} \in T_3$.\qed\\\\
{\bf Remark:} The results in Theorems 4 and 5 are true under the assumptions that $N_{\mathbb{Q}(\eta)/\mathbb{Q}}(\eta) \equiv 1\pmod 2$ i.e. constant term in the minimal polynomial of $\eta$ is an odd number and 
for cubes $3 \nmid \sigma$  and $3\nmid N_{\mathbb{Q}(\eta)/\mathbb{Q}}(\eta) $ . The problem still remains open without these assumptions.	

\end{document}